\documentclass[10p]{elsarticle}
\usepackage{amssymb, color}
\usepackage{amsmath}
 \usepackage{pifont}
 \usepackage[utf8]{inputenc}
\usepackage{setspace}
\journal{Journal of Dynamical and Control Systems}
\makeatletter
\def\LaTeX{\leavevmode L\raise.42ex
    \hbox{\kern-.3em\size{\sf@size}{0pt}\selectfont A}\kern-.15em\TeX}
\makeatother

\newcommand{\BibTeX}{{\rm B\kern-.05em{\sc
          i\kern-.025emb}\kern-.08em\TeX}}

\makeatletter
\def\@currentlabel{2.1}\label{e:dispaa}
\def\@currentlabel{2.21}\label{e:dispau}
\def\@currentlabel{2.22}\label{e:dispav}
\def\@currentlabel{2.23}\label{e:dispaw}
\def\@currentlabel{2.24}\label{e:dispax}
\def\theequation{\thesection.\@arabic\c@equation}
\makeatother

\renewcommand{\theequation}{\arabic{section}.\arabic{equation}}

\newcommand{\D }{\Delta }

\renewcommand{\l }{\lambda }

\newcommand{\var }{\varphi }

\newcommand{\be}{\begin{equation}}
\newcommand{\ee}{\end{equation}}
\newcommand{\ba}{\begin{align*}}
\newcommand{\ea}{\end{align*}}
\newcommand{\R}{\mathbb R}

\newcommand{\N}{\mathbb N}

\def \D{\Delta}

\newtheorem{thm}{Theorem} [section]
\newtheorem{lem}{Lemma} [section]
\newtheorem{prop}{Proposition} [section]

\newtheorem{defn}{Definition} [section]

\newtheorem{re}{Remark}[section]

\renewcommand{\theequation}{\thesection.\arabic{equation}}
\renewcommand{\thesection}{\arabic{section}}
\renewcommand{\theequation}{\thesection.\arabic{equation}}
\let\ssection=\section\renewcommand{\section}{\setcounter{equation}{0}\ssection}

\begin{document}
\begin{frontmatter}

\title{On the Classification of Stable Solutions of some elliptic equations in half-space.}
\author[mf1,mf2]{Foued Mtiri}
\ead{mtirifoued@yahoo.fr}
\author[ah1,ah2]{Abdelbaki Selmi\corref{cor1}}
\ead{Abdelbaki.Selmi@fsb.rnu.tn}
\author[mf3]{ Cherif Zaidi}
\ead{zaidicherif10@gmail.com}
\address[mf1]{Mathematics Department, Faculty of Sciences and Arts, King Khalid University, Muhayil Asir, Saudi Arabia.}
\address[mf2]{ Faculty of Sciences of Tunis, Elmanar University, Tunisia.}
\address[ah1]{ Mathematics Department, Northern Border University, Arar, Saudi Arabia.}
\address[ah2]{Universit\'{e} de Tunis, D\'{e}partement de Math\'{e}matiques, Facult\'{e} des Sciences de Bizerte, Zarzouna, 7021 Bizerte, Tunisia.}
\address[mf3]{Institut Sup\'{e}rieur des Sciences Appliqu\'{e}es et de Technologie de Kairouan, Avenue Beit El Hikma, 3100 Kairouan - Tunisia.}

\begin{abstract}
\noindent In this paper, we are concerned with stable solutions , possibly unbounded and sign-changing, of some semi-linear elliptic problem with mixed nonlinear boundary conditions. We establish the nonexistence of stable solutions, the main methods used are the Pohozaev identity, monotonicity formula of solutions together with a blowing down sequence.
\end{abstract}
\end{frontmatter}

 \section{Introduction and main results}
\setcounter{equation}{0}
This paper is devoted to the study of the following semi-linear elliptic problem
\begin{equation*}
(P)\ \  -\Delta u+\l u=|u|^{p-1}u\quad\mbox{in}\
\R^N_+,
\end{equation*}
where $\R^N_+=\{x=(x',x_N),\, x'\in\R^{N-1},\, x_N>0\}$, $N\geq2$, $p>1$ and $\l$ is a positive real parameter.
The motivation of studying such an equation is originated from the classical Lane-Emden equation
\begin{equation}\label{gs1}
-\D u=|u|^{p-1}u\,\,\mbox{in}\,\,\R^N
\end{equation}
A fundamental result on equation \eqref{gs1} is the celebrated Liouville-type theorem due to Gidas and Spruck \cite{GS}. The equation \eqref{gs1} has no positive classical solution if $1<p<\frac{N+2}{N-2}.$
 Moreover, it was also proved that the exponent $\frac{N+2}{N-2}$ is optimal, in the sense that problem \eqref{gs1} admits a positive solution for $p\geq\frac{N+2}{N-2}$ and $N\geq3$. Soon afterward, similar results were established in \cite{GS1} for positive solutions of the subcritical problem \eqref{gs1} in the upper half-space $\R^N_+$.
These results received wide attention as regards the theory itself and its applications. Particularly, when variational methods cannot be employed, one use them to establish a prior bound of solutions for general operator, and therefore existence of solutions may be dealt with via topological methods; see for instance \cite{GS, GS1}.
 \medskip
 
On the other hand, the idea of using the Morse index of a solution for a semilinear elliptic equation was first explored by Bahri and Lions \cite{BL} to get further proved that when $1<p<\frac{N+2}{N-2}$, no sign-changing solution  exists for \eqref{gs1}. To prove this result, they first deduced some integrable conditions on the solution based on finite Morse index; then they used the Pohozaev identity to prove the nonexistence result.
So, motivated by \cite{GS1}, they used blow-up argument to obtain a relevant $L^\infty$-bound for solutions of semilinear boundary value problems in bounded domain from the boundedness of Morse index (see also \cite{BL, Fa, HRS1, HHF, hyf}). We mention also that when the Palais-Smale; or the Cerami compactness conditions for the energy functional do not seem to follow readily, the proof of existence of solutions is essentially reduced to deriving  $L^\infty$-estimate from Liouville-type theorems via Morse index (see for instance \cite{DY, RR, RTT}). After these works, many authors investigated various Liouville type theorems for solutions with finite Morse indices in subcritical case such as problems with Neumann boundary condition, Dirichlet-Neumann mixed boundary and nonlinear boundary conditions (see \cite{MHA, HARS, HRS1, HRS2, WZ, XYU1, XYU2}). In the supercritical case, the finite Morse index solutions to the corresponding nonlinear problem \eqref{gs1} have been completely classified by Farina \cite{Fa}.
 \medskip
 
 A natural question is to understand more about finite Morse index solutions of the problem $(P)$ when $\l>0$. The case of Dirichlet boundary condition was studied by A. Selmi, A. Harrabi and C. Zaidi  in \cite{HSZ1}. They prove various Liouville type theorems for stable solutions possibly unbounded and sign-changing. While when the half space is replaced by strip domain similar result was obtained in \cite{HSZ2}. 
 \medskip
 
 In this paper we study solutions, possibly unbounded and sign-changing of the following mixed problems with mixed boundary value conditions
\begin{equation}\label{11}
 \left\{\begin{array}{lllll} -\Delta u+\l u=|u|^{p-1}u,\quad \quad\mbox{in}\ 
\R^N_+,\\\\
\frac{\partial u}{\partial \nu}=|u|^{q-1}u, \quad \quad  \quad \quad  \quad \mbox{on}\ 
\Sigma_1,\\\\
 \frac{\partial u}{\partial \nu}=0, \quad \quad \quad \quad  \quad \quad  \quad \quad \mbox{on}\ 
\Sigma_2,
\end{array}
\right.
\end{equation}
or
\begin{equation}\label{12}
\left\{\begin{array}{lllll} -\Delta u+\l u=|u|^{p-1}u \quad \quad\mbox{in}\
\R^N_+,\\\\
\frac{\partial u}{\partial \nu}=|u|^{q-1}u \quad \quad  \quad \quad  \quad\mbox{on}\
\Sigma_1,\\\\
 u=0 \quad \quad \quad \quad  \quad \quad  \quad\quad \mbox{on}\
\Sigma_2,
\end{array}
\right.
\end{equation}
where
 \begin{align*}
 &\R^N_+=\{x=(x',x_N),\, x'\in\R^{N-1},\, x_N>0\},\\
 &\Sigma_1=\{x=(x_1,\ldots x_N)\in\R^N; x_N=0, x_1>0 \},\\
 &\Sigma_2=\{x=(x_1,\ldots x_N)\in\R^N; x_N=0, x_1<0 \}.
\end{align*}
In the case $\l=0$, Liouville-type theorems and properties for problem \eqref{11} or \eqref{12} has been studied by \cite{WZ, XYU2}. In \cite{WZ} it was shown that when  $\l=0$ there is no nontrivial bounded solution of problem \eqref{11} or \eqref{12} with finite morse index, provided that
$$1<p\leq\frac{N+2}{N-2}, 1<q\leq\frac{N}{N-2}\,\, \mbox{and}\,\,(p,q)\ne(\frac{N+2}{N-2},\frac{N}{N-2}).$$
Recently in \cite{HRa} Harrabi and Rahal have improved the results in \cite {XYU2} for a large class of the exponents $p$ and $q$.
\medskip

Our aim in this work is to analyze the influence of the linear term $\l u$ when $\l >0$ in order to classify regular stable solutions possibly unbounded and sign-changing. We will prove various Liouville type theorems for $C^2$ solutions which are stable or stable outside a compact set. our analysis reveals the nonexistence of nontrivial finite Morse index solutions for all $p>1$. Following \cite{DDWW, DH, HR, HSZ1, FP}, we establish a variant of the monotonicity formula to provide nonexistence results in the supercritical range. However, in the subcritical case, we require a restrictive condition on $u$ namely $$|u|_{L^\infty(\R^N_+)}^{p-1}< \l(\frac{N}{2}-\frac{N-1}{q+1})(\frac{N}{p+1}-\frac{N-1}{q+1})^{-1},$$ which allows-us to provide nonexistence result from the Pohozaev identity.
\medskip

We mention that Liouville type theorems in unbounded domains play a
crucial role to obtain a priori $L^\infty$-bounds for solutions of semilinear boundary value problems in bounded
domain, see \cite{GS} for the case of positive solutions and \cite{BL}, \cite{Da2}, \cite{Fa} and \cite{HRS1} for sign changing solutions and having finite
Morse index. Also, nonlinear Liouville type theorems combined with the degree type arguments, are useful
to obtain the existence of solutions (see for instance \cite{BDN}).
\medskip

In order to state our results we need to recall the following.
\begin{defn}
We say that a solution $u$ of \eqref{11} belonging to $C^2(\overline{\R^N_+})$
\medskip

$\bullet$ is stable if
\begin{equation}\label{st}
Q_u(\varphi)=\int_{\R^N_+}|\nabla \varphi|^2-q\int_{\Sigma_1}|u|^{q-1}\var^2+\l\int_{\R^N_+}\varphi^2-p\int_{\R^N_+}|u|^{p-1}\var^2\geq 0\,\,\,\forall\varphi\in C^1_c(\overline{\R^N_+}),
\end{equation}
$\bullet$ has a Morse index equal to $K\geq 1$ if $K$ is the maximal dimension of a subspace $X_K$ of $C^1_c(\overline{\R^N_+})$ such that $Q_u(\varphi)<0$ for any $\varphi\in X_K\backslash\{0\}$,
\medskip

$\bullet$ is stable outside a compact set $\mathcal{K}$ if $Q_u(\varphi)\geq0$ for any $\varphi\in C^1_c(\overline{\R^N_+}\backslash\mathcal{K})$.
\end{defn}
Similarly, we say that a solution $u$ of \eqref{12} belonging to $C^2(\overline{\R^N_+})$ is stable
(respectively, stable outside a compact set $\mathcal{K}$) if $Q_u(\varphi)\geq0$ for every $\varphi\in C^1_c(\R^N_+\cup\Sigma_1)$ (respectively, $\varphi\in C^1_c\big((\R^N_+\cup\Sigma_1)\backslash \mathcal{K}$\big).

\begin{re}
\begin{itemize}
\item[ (i).] Clearly a solution is stable if and only if its Morse index is equal to zero.
\item[ (ii).] Any finite Morse index solution $u$ is stable outside a compact set $\mathcal{K}\subset\R^N_+$. Indeed there exist $K\geq1$ and $X_K:=$ span$\{\varphi_1,\ldots,\varphi_K\}\subset C^1_c(\R^N_+)$ such that $Q_u(\varphi)<0$ for any $\varphi\in X_K\backslash\{0\}$. Then, $Q_u(\varphi)\geq0$ for every $\varphi\in C^1_c(\overline{\R^N_+}\backslash\mathcal{K})$, where $\mathcal{K}:=\cup^K_{j=1}\mbox{supp}(\varphi_j)$.
\end{itemize}
\end{re}
The first result of this paper is
\begin{thm}\label{t1}
Let $u\in C^2(\overline{\R^N_+})$, be a stable solution of \eqref{11} or \eqref{12} and  $\l>0$. Assume that
\begin{itemize}
\item[ 1.] $p>1$, $q>1$ if $N=2$,
\item[ 2.] $p\in (1,\frac{N+2}{N-2}]$ and $q>1,$ or $p>\frac{N+2}{N-2}$ and $2q-p-1\geq0$ if $N\geq3$.
\end{itemize}
Then $u\equiv 0$.
\end{thm}
\begin{re}
In statement (2) of Theorem \ref{t1} if $p>\frac{N+2}{N-2}$ then $2q-p-1\geq0$ is satisfied if $p\leq q$. If $p> q$ the condition $2q-p-1\geq0$ implies that necessarily $q> \frac{N}{N-2}$.
\end{re}
In the case of solutions of \eqref{11} or \eqref{12} which are stable outside a compact set of $\R^N_+$ we prove that
\begin{thm}\label{t2}
Let $N\geq3$, $\l>0$ and $u\in C^2(\overline{\R^N_+})$ be a solution of \eqref{11} or \eqref{12} which is stable outside a compact set.  Assume that $$|u|_{L^\infty(\R^N_+)}^{p-1}<\l (\frac{N}{2}-\frac{N-1}{q+1})(\frac{N}{p+1}-\frac{N-1}{q+1})^{-1},\quad p\in(1,\frac{N+2}{N-2})\quad \mbox{and} \;\;\; q> \frac{N}{N-2},$$
 then  $u\equiv0.$
\end{thm}
\begin{thm}\label{t3}
Let $N\geq3$, $\l>0$ and $u\in C^2(\overline{\R^N_+})$, be a solution of \eqref{11} or \eqref{12} which is stable outside a compact set of $\R^N_+$. Assume that $p\geq\frac{N+2}{N-2}$ with $2q-p-1\geq0$, then $u \equiv0$.
\end{thm}
The proof of Theorem \ref{t2} or \ref{t1} uses a version of monotonicity formula of equation \eqref{11} or \eqref{12}. We mention that the monotonicity formula is a powerful tool to understand supercritical elliptic equations or systems. This approach has been used successfully for the Lane-Emden equation in \cite{FP}.
\medskip

For $R>0$, denote $A_R=\{ R<|x|<2 R\}$, $B_R^+=B_R\cap \R^N_+$ and $\partial B_R^+=\partial B_R\cap\R^N_+$ where $B_R$ is the open ball centered at the origin and with radius $R$. The key step of the proofs of Theorems \ref{t1}, \ref{t2} and \ref{t3} is the following integral estimates which are useful in the subcritical, critical and supercritical cases.
\begin{prop}\label{pp}
Let $u\in C^2(\overline{\R^N_+})$ a solution of \eqref{11} or \eqref{12} which is stable outside a compact set $\mathcal{K} \subset  B^+_{R_0} $ for some $R_0>0$. Then for all $R> 2R_0,$ we have.
\begin{eqnarray}\label{int1}
\int_{B^+_{ R }} |\nabla u|^2+\int_{ B^+_{ R }} \l u^2+ \int_{ B^+_{ R }}|u|^{p+1}+ \int_{ \Sigma_1\cap B^+_{ R }}|u|^{q+1} \leq C_0+ CR^{-2 }\int_{ A_R} u^{2},
\end{eqnarray}
 and
 \begin{eqnarray}\label{int11}
\int_{B^+_{R} } |\nabla u|^2+\int_{ B^+_{ R }} \l u^2+ \int_{B^+_{R} }|u|^{p+1}+ \int_{ \Sigma_1\cap B^+_{ R }}|u|^{q+1} \leq C_0+ CR^{N-2\frac{p+1}{p-1}}.
\end{eqnarray}
Here $C=C(N,p)$ and $C_0=C_0(u,R_0,N,p)$ are positive constants independent of $R$. Furthermore, if $u$ is a stable solution, then \eqref{int1} and \eqref{int11} hold with $C_0=0$.
\end{prop}
\begin{re}\label{bou-sta}
Observe that classification of {\bf bounded stable} solutions of \eqref{11} or \eqref{12} follows immediately from \eqref{int1} for all $p>1$. In fact, since $u$ is bounded, then from \cite{RR} we can find a positive constant $C$ and a sequence $R_n \to \infty \mbox{ as } n\to \infty$  such that
\begin{eqnarray}\label{RR}
 \int_{\R^N_+\cap B_{2R_n}} u^{2} \leq C \int_{\R^N_+\cap B_{R_n}} u^{2}.
\end{eqnarray}
 According to \eqref{int1} (with $C_0=0$), we derive
\begin{eqnarray*}
 \l  \int_{\R ^N_+\cap B_{ R_n }}u^{2}  \leq  C R_n^{-2 } \int_{\R^N_+\cap B_{R_n}} u^{2}  .
\end{eqnarray*}
As $R_n \to \infty \mbox{ as } n\to \infty$ , there exists $n_0\in \N$ such that $ C R_n^{-2 } <\frac {\l}2$ for all $n>n_0$, and therefore we deduce that
\begin{eqnarray*}
   \int_{\R^N_+\cap B_{ R_n }}u^{2}  \leq   0,\; \forall n>n_0,
\end{eqnarray*}
which implies that $u\equiv 0$. Point out that for $\l>0$, we cannot always use the doubling lemma technique as in the case $\l=0$ (see \cite{ PP}), to reduce the classification for only bounded solutions of \eqref{11} or \eqref{12}. Therefore, we shall pay special attention to the delicate case of unbounded solutions by exploiting the following variant of the monotonicity formula.
\end{re}
For $\tau>0$, define the function $u^\tau$ by
$$u^\tau(x)=\tau^\frac{2}{p-1}u(\tau x), \;\;\mbox{for}\;\; x\in \R^N_+,$$
then we have
\begin{prop}\label{po}
Let $u$ $\in C^2(\overline{\R^N_+})$ be a solution of \eqref{11} or \eqref{12} and $\tau>0$ be a constant. Set
\begin{align}\label{M1}
\begin{split}
E(u,\tau)&=\int_{B_{1}\cap\R^N_+}\big(\frac{1}{2}|\nabla u^\tau|^{2}+\frac{\lambda\tau^2}{2}|u^\tau|^{2}-\frac{1}{p+1}|u^\tau|^{p+1}\big)\\
&-\frac{\tau^{1-2\frac{q-1}{p-1}}}{q+1}\int_{\Sigma_1\cap B_1}|u^\tau|^{q+1}+\frac{1}{p-1}\int_{\partial B_{1}\cap\R^N_+}|u^\tau|^{2}.
\end{split}
\end{align}
Then
\begin{align}\label{M2}
\frac{d E}{d\tau}= \tau \int_{\partial B_{1}\cap\R^N_+}(\frac{d u^\tau}{d\tau})^2ds+\lambda\tau \int_{B_{1}\cap\R^N_+}(u^\tau)^{2}dx+\frac{2q-p-1}{(p-1)(q+1)}\tau^{-2\frac{q-1}{p-1}}\int_{\Sigma_1\cap B_1}|u^\tau|^{q+1}
\end{align}
for all $p,q>1$. Furthermore, $E$ is a nondecreasing function of $\tau$ if $2q-p-1\geq0$ .
\end{prop}
\medskip

This paper is organized as follows. In section $2$ we prove Propositions \ref{pp} and \ref{po}. The proofs of Theorems \ref{t1} and \ref{t2}-\ref{t3} are given respectively in sections $3$ and $4.$

\medskip

\section{Proofs of Propositions  \ref{pp} and \ref{po}}
\setcounter{equation}{0}
\medskip

\textbf{ Proof of Proposition \ref{pp}}. Let $u\in C^2(\overline{\R^N_+})$  be a solution of \eqref{11} or \eqref{12} which is stable outside a compact set $\mathcal{K} \subset B^+_{R_0}$. For all $R> 2R_0 $, define the family of test functions $\psi= \psi _{(R,R_0) }\in C^{1}_c(\R^n)$ satisfying
\begin{eqnarray} \label{eq1} \begin{cases}
0 \leq \psi \leq 1,\\
\psi\equiv 1 \;\;\text{if \; $2R_0<|x|< R$},\; \psi\equiv 0\;\,\text{if\; $|x|<R_0$  \;\mbox{or }\; $|x|> 2 R$,}\\
|\nabla \psi|\leq C R^{-1}\;\text{if\; $ R<|x|<2R$}.
\end{cases}
\end{eqnarray}
If $u$ is a stable solution, we use $\psi=\psi_{R}$ verifying \eqref{eq1} with $R_0=0$ (that is $\psi=1$ if  $|x|< R$). For $m > 1$
we have $u\psi^m \in C_c^1( \overline{B^+_{2R}\backslash B^+_{R_0}})$, then it can be chosen as a test function and in view of \eqref{st}, we have 
\begin{eqnarray}\label{e1}
\displaystyle\int_{B^+_{2R}}|\nabla (u\psi^m)|^2+\l\int_{B^+_{2R}} u^{2}\psi^{2m} dx -p\int_{B^+_{2R}}  |u|^{p+1}\psi^{2m} dx -q \int_{\Sigma_1\cap B^+_{2R}}  |u|^{q+1}\psi^{2m} \geq 0, \;\; \forall\;\; R>2R_0.
\end{eqnarray}
Multiply \eqref{11} or \eqref{12} by $u\psi^{2m}$ and integrate by parts and taking into account that $$\nabla u.\nabla (u\psi^{2m})=|\nabla (u\psi^m)|^2-u^2|\nabla (\psi^m)|^2 ,$$
so, we get 
\begin{align}
\begin{split}
\label{e11}
& \;\int_{B^+_{2R}}|\nabla (u\psi^m)|^2 +\l\int_{B^+_{2R}} u^{2}\psi^{2m} dx
- \int_{B^+_{2R}} |u|^{p+1}\psi^{2m} dx\\
& \;
=- \int_{\Sigma_1\cap B^+_{2R}} |u|^{q+1}\psi^{2m} dx+\int_{B^+_{2R}}u^2|\nabla (\psi^m)|^2 .
\end{split}
\end{align}
Combine \eqref{e1} and \eqref{e11}, then for all  $R>2R_0$ we derive
\begin{align}
\begin{split}\label{e111}
& \;\displaystyle\int_{B^+_{2R}}|\nabla (u\psi^m)|^2 +\l\int_{B^+_{2R}} u^{2}\psi^{2m} dx+\int_{B^+_{2R}}  |u|^{p+1}\psi^{2m} dx
+\int_{\Sigma_1\cap B^+_{2R}}  |u|^{q+1}\psi^{2m} dx \\
& \;\leq C_p\displaystyle\int_{B^+_{2R}}u^2|\nabla (\psi^m)|^2,
\end{split}
\end{align}
which implies that
\begin{align}
\begin{split}\label{zz}
& \;\displaystyle\int_{B^+_{2R}}|\nabla (u\psi^m)|^2 +\l\int_{B^+_{2R}} u^{2}\psi^{2m} dx+ \int_{B^+_{2R}}  |u|^{p+1}\psi^{2m} dx
 +\int_{\Sigma_1\cap B^+_{2R}}  |u|^{q+1}\psi^{2m} dx\\
& \;\leq C_0+ C_{p,m}R^{-2}\displaystyle\int_{A^+_{R}}u^2 \psi^{2m-2},
\end{split}
\end{align}
where $C_0=0$ if $u$ is a stable solution. Hence inequality \eqref{int1} follows from the last inequality. Using now Young's inequality and choose $m=\frac{p+1}{p-1}> 1 $, we deduce that
$$ R^{-2}\displaystyle\int_{A^+_{R}}u^2\psi^{2m-2} \leq R^{-2}\displaystyle\int_{B^+_{2R}}u^2\psi^{2m-2} \leq C_{p,N}R^{N-2\frac{p+1}{p-1}}+\frac{2}{p+1}\displaystyle\int_{B^+_{2R}}u^{p+1}\psi^{2m}.$$
 Insert the last inequality into the right hand-side of \eqref{zz}, so derive the main integral estimate \eqref{int11}. The proof of Proposition \ref {pp} is thereby completed. \qed
 \medskip
 
 \textbf{ \textbf{ Proof of Proposition \ref{po}}}.
Since $u$ is a solution of \eqref{11}, it follows that $u^\tau$ satisfies
\begin{equation}\label{mp}
\quad\left\{\begin{array}{lllllllll} -\Delta u^\tau+\l \tau^2 u^\tau=|u^\tau|^{p-1}u^\tau \quad\quad\quad\mbox{in}\
\R^N_+,\\\\
\frac{\partial u^\tau}{\partial \nu}=\tau^{1-2\frac{q-1}{p-1}}|u^\tau|^{q-1}u^\tau \quad\quad\quad\quad \mbox{on}\ \ \Sigma_1,\\\\
 \frac{\partial u^\tau}{\partial \nu}=0 \quad\quad\quad\quad\quad\quad\quad\quad\quad\quad\quad\mbox{on}\ \Sigma_2.
\end{array}
\right.
\end{equation}

Let
\begin{equation}\label{mj}
\overline{E}(u,\tau)=\int_{B_{1}\cap\R^N_+}\bigg(\frac{1}{2}|\nabla u^\tau|^{2}+\frac{\l}{2}\tau^2(u^\tau)^2-\frac{1}{p+1}|u^\tau|^{p+1}\bigg),
\end{equation}
then
\begin{equation}\label{mc}
\frac{d}{d\tau}\overline{E}(u,\tau)=\int_{B_{1}\cap\R^N_+}\nabla u^\tau\nabla\frac{du^\tau}{d\tau}+\l\tau\int_{B_1}(u^\tau)^2+\l\tau^2\int_{B_1}u^\tau\frac{du^\tau}{d\tau}-\int_{B_1}|u^\tau|^{p-1}u^\tau\frac{du^\tau}{d\tau}.
\end{equation}
Integrating by parts and using the fact that $u^\tau$ is a solution of \eqref{mp}, we get
\begin{equation}\label{ml}
\frac{d}{d\tau}\overline{E}(u,\tau)=\int_{\partial B_{1}\cap\R^N_+}\frac{\partial u^\tau}{\partial r}\frac{du^\tau}{d\tau}+\l\tau\int_{B_1}(u^\tau)^2+\tau^{{1-2\frac{q-1}{p-1}}}\int_{\Sigma_1\cap B_1}|u^\tau|^{q-1}u^\tau\frac{du^\tau}{d\tau} .
\end{equation}
In what follows, we express all derivatives of $u^\tau$ in the $r=|x|$ variable in terms of derivatives in the $\tau$ variable. In the definition of $u^\tau$, directly differentiating in $\tau$ gives
\begin{equation}\label{mn}
\tau\frac{du^\tau}{d\tau}=\frac{2}{p-1} u^\tau + r \frac{\partial u^\tau}{\partial r}.
\end{equation}
From \eqref{ml} and \eqref{mn}, we obtain
\begin{align}
\begin{split}
\label{ms}
\frac{d}{d\tau}\overline{E}(u,\tau)&=\tau \int_{\partial B_{1}\cap\R^N_+}(\frac{du^\tau}{d\tau})^2-\frac{1}{p-1}\int_{\partial B_{1}\cap\R^N_+}\frac{d(u^\tau)^2}{d\tau}
\\&+\lambda\tau\int_{B_{1}\cap\R^N_+}(u^\tau)^2+\tau^{{1-2\frac{q-1}{p-1}}}\int_{\Sigma_1\cap B_1}|u^\tau|^{q-1}u^\tau\frac{du^\tau}{d\tau}.
\end{split}
\end{align}
On the other hand we have
\begin{align}
\begin{split}
\label{mv}
 & \;\tau^{{1-2\frac{q-1}{p-1}}}\int_{\Sigma_1\cap B_1}|u^\tau|^{q-1}u^\tau\frac{du^\tau}{d\tau}\\
=& \;\frac{d}{d\tau}\left[ \frac{\tau^{{1-2\frac{q-1}{p-1}}}}{q+1}\int_{\Sigma_1\cap B_1}|u^\tau|^{q+1}\right]-\frac{(p+1-2q)\tau^{-2\frac{q-1}{p-1}}}{(p-1)(q+1)}\int_{\Sigma_1\cap B_1}|u^\tau|^{q+1} .
\end{split}
\end{align}
Combining \eqref{mj}, \eqref{ms} and \eqref{mv}, we get \eqref{M1} and \eqref{M2}.
\smallskip

Concerning problem \eqref{12} the proof can be obtained with only minor modifications. Since $u$ is a solution of \eqref{12}, then $u^\tau$ satisfies
\begin{equation}\label{mp1}
\quad\left\{\begin{array}{llllll} -\Delta u^\tau+\l \tau^2 u^\tau=|u^\tau|^{p-1}u^\tau,\quad\quad\mbox{in}\
\R^N_+,\\
\frac{\partial u^\tau}{\partial \nu}=\tau^{1-2\frac{q-1}{p-1}}|u^\tau|^{q-1}u^\tau, \quad \quad \quad \quad   \mbox{on}\ \ \Sigma_1,\\
 u^\tau=0, \quad \quad \quad \quad  \quad \quad  \quad \quad\mbox{on}\ \Sigma_2.
\end{array}
\right.
\end{equation}
Integrating by parts in \eqref{mc} and using the fact that $u^\tau$ is a solution of \eqref{mp1}, we get
\begin{align}
\begin{split}
\label{mz}
\frac{d}{d\tau}\overline{E}(u,\tau)&=\int_{\partial B_{1}\cap\R^N_+}\frac{\partial u^\tau}{\partial r}\frac{du^\tau}{d\tau}+\l\tau\int_{B_1}(u^\tau)^2+\int_{(\Sigma_1\cup\Sigma_2)\cap B_1}\frac{\partial u^\tau}{\partial\nu}\frac{du^\tau}{d\tau}\\
&=\int_{\partial B_{1}\cap\R^N_+}\frac{\partial u^\tau}{\partial r}\frac{du^\tau}{d\tau}+\l\tau\int_{B_1}(u^\tau)^2+\tau^{{1-2\frac{q-1}{p-1}}}\int_{\Sigma_1\cap B_1}|u^\tau|^{q-1}u^\tau\frac{du^\tau}{d\tau}.
\end{split}
\end{align}
Indeed, we have $u^\tau\equiv0$ in $\Sigma_2\cap B_1$ for all $\tau>0$, then $\frac{du^\tau}{d\tau}=0$ in $\Sigma_2\cap B_1$. The rest of the prof is unchanged.
Now, since $\lambda$ is a positive scalar and $2q-p-1\geq0$ we have that $E$ is a nondecreasing function of $\tau$. This completes the proof of Proposition \ref{po}. \qed

\section{Proof of Theorem \ref{t1}}
\setcounter{equation}{0}
\medskip

Let $u\in C^2(\overline{\R^N_+})$ be a stable solution of \eqref{11} or \eqref{12}.
If $N=2$ from \eqref{int11} we have
$$\int_{B_R^+}|u|^{p+1}\leq CR^{-\frac{4}{p-1}},$$
which yields that $u\equiv0$ for all $p>1$. 

\medskip

If $N\geq3,$ three cases may occur.
\medskip

 {\bf Case $1.$ If $1<p<\frac{N+2}{N-2}$}. From the main integral estimate \eqref{int11} of Proposition \ref{pp} (with $C_0=0$), we have
\begin{equation*}
\int_{B_R^+}|u|^{p+1}\leq CR^{N-2\frac{p+1}{p-1}}, \; \forall R>0.
\end{equation*}
As $N-2\frac{p+1}{p-1}<0$ if $1<p<\frac{N+2}{N-2}$, we may readily see that $u\equiv 0$.
\medskip

{\bf Case $2.$ If $p=\frac{N+2}{N-2}$}. Apply again \eqref{int11}, we derive
\begin{equation*}
\int_{\R_+^N}|u|^{\frac{ 2N}{N-2}} < \infty \mbox{ and } \int_{A_R\cap \R_+^N}|u|^{\frac{ 2N}{N-2}}\to 0 \mbox{ as } R \to \infty,
\end{equation*}
where $A_R=\{ R<|x|< 2R \}$. Invoking now \eqref{int1} of Proposition \ref{pp} (with $C_0=0$) and applying  H\"{o}lder's inequality, we derive that
\begin{equation*}
\int_{B_R^+}|u|^{\frac{ 2N}{N-2}}\leq CR^{-2}\int_{A_R\cap \R_+^N }|u|^{2} \leq C\left(\int_{A_R\cap \R_+^N}|u|^{\frac{ 2N}{N-2}}\right)^\frac{ N-2}{N}, \; \forall R>0.
\end{equation*}
Hence, we deduce that $\displaystyle\int_{\R_+^N}|u|^{\frac{ 2N}{N-2}}=0$ and so $u\equiv 0$.
\medskip

{\bf Case $3$. If $p>\frac{N+2}{N-2}$}. This case needs more involving analysis using the powerful monotonicity formula. The fact that $\lambda >0$ allows-us to provide nonexistence result for all supercritical exponent $p$.
\smallskip

Substitute $R$ by $\tau R$ in \eqref{int11} (with $C_0=0$), then $\forall R>0 \mbox{ and }\tau >1$, one has
\begin{equation}\label{sca}
\int_{B_{\tau R}^+}\bigg(|\nabla u|^2+|u|^{p+1}+\l|u|^{2}\bigg)+\int_{\Sigma_1\cap B_{\tau R}^+}|u|^{q+1}\leq C\tau^{N-2\frac{p+1}{p-1}}R^{N-2\frac{p+1}{p-1}} .
\end{equation}
Scaling back\footnote{ In fact, multiply \eqref{sca} by $\tau^{2\frac{p+1}{p-1}}$ and use the change of variable $x=\tau y$.}, we obtain
\begin{equation}\label{rr}
\int_{B_R^+}\bigg(|\nabla u^\tau|^2+|u^\tau|^{p+1}+\l\tau^2|u^\tau|^{2}\bigg)+\tau^{1-2\frac{q-1}{p-1}}\int_{\Sigma_1\cap B_R^+}|u^\tau|^{q+1}\leq CR^{N-2\frac{p+1}{p-1}}, \; \forall R>0, \forall \tau >1.
\end{equation}
So, $(u^\tau)_{\tau \geq 1}$ is uniformly bounded in $H^1\cap L^{p+1}(B_R^+)$ and $(\tau^{\frac{p+1-2q}{(p-1)(q+1)}}u^\tau)_{\tau \geq 1}$ is uniformly bounded in $L^{q+1}(\Sigma_1\cap B_R)$  for any $R>0$. Consequently, from a standard diagonal subsequence argument, we can find a sequence $(u^{\tau_j})$ which converges strongly in $L^{2}(B_R^+)$ and weakly in $H^1\cap L^{p+1}(B_R^+)$ to some function $u^\infty$ for every $R>0$, as $\tau_j\rightarrow +\infty$. Multiply the following equation
\begin{eqnarray*}
  u^\tau=\lambda^{-1}\tau^{-2}(\D u^\tau + |u^\tau|^{p-1}u^\tau) \; &\text{ in \; $\R^{N}_{+}$},
\end{eqnarray*}
by $\phi \in C_c^1(\overline{B^+_R})$, then integrate by parts and using H\"{o}lder's inequality, we derive from \eqref{rr} that
$$\bigg|\int_{B_R^+}u^\tau \phi\bigg|\leq \dfrac{C_{R,\phi}}{\tau^{2}},  \mbox{ where } C_{R,\phi}> 0. $$
Consequently, $(u^{\tau_j})$  converges weakly to $0$ in $L^{2}(B_R^+)$, so $u^\infty=0$ and then $(u^{\tau_j})$ converges strongly to $0$ in $L^{2}(B_R^+)$.
Invoking now inequality \eqref{int1} of Proposition \ref{pp} where we substitute $R$ by $\tau R$ we obtain
\begin{equation*}
\int_{B_{\tau R}^+}(|\nabla u|^2+|u|^{p+1}+\l u^{2})+\int_{ \Sigma_1\cap B^+_{ \tau R }}|u|^{q+1}\leq CR^{-2}\tau^{-2}\int_{B_{2\tau R}^+}u^2.
\end{equation*}
Thus, scaling back, we deduce
\begin{equation}\label{lp}
\int_{B_{ R}^+}(|\nabla u^\tau|^2+|u^\tau|^{p+1}+\l\tau^2|u^\tau|^{2})+\tau^{1-2\frac{q-1}{p-1}}\int_{\Sigma_1\cap B^+_R}|u^\tau|^{q+1}\leq CR^{-2}\int_{B_{2 R}^+}(u^\tau)^2.
\end{equation}
As a consequence, $(u^{\tau})$ converges strongly to $0$ in $H^1\cap L^{p+1}(B_R^+)$ and $(\tau u^{\tau})$ converges strongly to $0$ in $L^{2}(B_R^+)$ as $\tau\to \infty$. Moreover, we have
\begin{equation}\label{lb}
\lim_{\tau\rightarrow +\infty}E_2(u,\tau)=0,
\end{equation}
where $$E_2(u,\tau)=\int_{B_{1}^+}\bigg(\frac{1}{2}|\nabla u^\tau|^{2}+\frac{\lambda\tau^2}{2}|u^\tau|^{2}-\frac{1}{p+1}|u^\tau|^{p+1}\bigg)-\frac{\tau^{1-2\frac{q-1}{p-1}}}{q+1}\int_{\Sigma_1\cap B_1}|u^\tau|^{q+1}.$$
We claim that the same holds true for $E$. To see this, simply observe that since $E$ is nondecreasing,
\begin{eqnarray*}
E(u,\tau)&\leq &\frac{1}{\tau}\int_{\tau}^{2\tau}E(u,t)dt\nonumber\\&=&\frac{1}{\tau}\int_{\tau}^{2\tau}E_2(u,t)dt+\frac{\tau^{-1}}{p-1}\int_{\tau}^{2\tau}\int_{\partial B_{1}^+}|u^t|^2\nonumber\\&\leq&\sup_{t\geq\tau}E_2(u,t)+C\int_{B_{2}^+} |u^\tau|^2.
\end{eqnarray*}
Thanks to this, we deduce from \eqref{lb} that
$$\lim_{\tau\rightarrow +\infty}E(u,\tau)=0.$$
In addition, since $u\in C^2(\overline{\R_+^N}),$ one easily verifies that $E(u,0)=0$.
As a consequence $E(u,\tau)\equiv 0$, and therefore $\frac{dE}{d\tau}=0$. Then from \eqref{M2}
$$ \int_{B_{1}^+}|u^\tau|^{2}dx=0,\, \forall \tau >1.$$
We readily deduce that $u\equiv 0$. The proof is completed. \qed

\section{Proof of Theorems \ref{t2}-\ref{t3}}
\setcounter{equation}{0}
\medskip

In this section, we will prove Theorems \ref{t2}-\ref{t3}. For this, we need the following well known Pohozaev identity.
\begin{lem}\label{l4}
Let $u$ $\in C^2(\overline{\R^N_+})$ be a solution of \eqref{11} or \eqref{12}. Then the following identity holds:
\begin{align}
\begin{split}\label{j}
& \;\frac{N-2}{2}\int_{B_{R}^+}|\nabla u|^{2}+\frac{N}{2}\l\int_{B_{R}^+}u^2-\frac{N}{p+1}\int_{B_{R}^+}|u|^{p+1} -\frac{N-1}{q+1}\int_{B_{R}\cap \Sigma_1}|u|^{q+1}dx'\\
=& \;\frac{R}2\int_{\partial B_{R}^+}|\nabla u|^2ds-R\int_{\partial B_{R}^+} \Big|\frac{ \partial u}{\partial\nu}\Big|^{2}ds +\frac{R}2\l\int_{\partial B_{R}^+}u^2 ds-\frac{R}{p+1}\int_{\partial B_{R}^+}| u|^{p+1}ds\\
& \;
-\frac{R}{q+1}\int_{\partial B_{R}\cap\Sigma_1}|u|^{q+1}ds.
\end{split}
\end{align}
\end{lem}

\noindent{\bf Proof.} The proof of this lemma is standard, we give it here for completeness. We deal only
with problem \eqref{11}. The proof for problem \eqref{12} is almost the same except that different
boundary value conditions were used. We omit the details.
\smallskip

Multiplying the equation \eqref{11} by $\langle x,\nabla u\rangle$ and integrating on $B_R^+$, then a direct computation shows that
\begin{align}\begin{split}\label{p1}
\int_{B_{R}^+}|u|^{p-1}u\langle x,\nabla u\rangle dx&=\frac{1}{p+1}\int_{B_{R}^+}\langle x,\nabla |u|^{p+1}\rangle dx\\
&\;=-\frac{N}{p+1}\int_{B_{R}^+}|u|^{p+1}dx+\frac{1}{p+1}\int_{\partial B_{R}^+}\langle x,\nu\rangle| u|^{p+1}ds\\
&\;=-\frac{N}{p+1}\int_{B_{R}^+}|u|^{p+1}dx+\frac{R}{p+1}\int_{\partial B_{R}^+}| u|^{p+1}ds.
\end{split}
\end{align}
Similarly, we have
\begin{align}\label{p3}
\int_{B_{R}^+}u\langle x,\nabla u\rangle dx
=-\frac{N}2\int_{B_{R}^+}u^2dx+\frac{R}2\int_{\partial B_{R}^+}u^2 ds.
\end{align}
Next, we deduce
\begin{align}\begin{split}\label{p2}
\int_{B_{R}^+}-\D u\langle x,\nabla u\rangle dx&=\int_{B_{R}^+}\nabla u\nabla (\langle x,\nabla u\rangle) dx-\int_{\partial B_{R}^+}\frac{\partial u} {\partial\nu}\langle x,\nabla u\rangle ds-\int_{B_{R}\cap\partial \R^N_+}\frac{\partial u} {\partial\nu}\langle x,\nabla u\rangle dx'\\
&=\int_{B_{R}^+}|\nabla u|^2dx+\frac{1}{2}\int_{B_{R}^+}\langle x,\nabla(|\nabla u|^2)\rangle dx-R\int_{\partial B_{R}^+} \Big|\frac{ \partial u}{\partial\nu}\Big|^{2}ds\\&\;\;\;-\frac{1}{q+1}\int_{B_{R}\cap\Sigma_1}\langle x',\nabla_{x'} (|u|^{q+1})\rangle dx'\\&= -\frac{N-2}{2}\int_{B_{R}^+}|\nabla u|^2dx+\frac{R}2\int_{\partial B_{R}^+}|\nabla u|^2ds-R\int_{\partial B_{R}^+} \Big|\frac{ \partial u}{\partial\nu}\Big|^{2}ds\\&\;\;\;+\frac{N-1}{q+1}\int_{B_{R}\cap \Sigma_1}|u|^{q+1}dx' -\frac{R}{q+1}\int_{\partial B_{R}\cap\Sigma_1}|u|^{q+1} ds.
\end{split}
\end{align}
Then \eqref{j} follows immediately from \eqref{p1}, \eqref{p3} and \eqref{p2}.\qed
\medskip

\subsection{\bf Proof of Theorem \ref{t2}.}
\medskip

We only prove the conclusion for problem \eqref{11}, the proof for problem \eqref{12} is the same.
Working by contradiction. Suppose that $u\neq0$ then, by virtue of the main integral estimate \eqref{int11} of Proposition \ref{pp} (with $C_0> 0$) we have
$$\int_{\Sigma_1\cap B_R^+}|u|^{q+1}, \; \int_{\R^N_+}|u|^{p+1}, \; \int_{\R^N_+}u^2 \mbox{ and } \int_{\R^N_+}|\nabla u|^2<\infty;\quad \mbox{ for any }\quad p\in (1,\frac{N+2}{N-2}).$$
Consequently, we can find a sequence $R_n \to \infty$  such that
\begin{align*}
&\frac{R_n}2\int_{\partial B_{Rn}^+}|\nabla u|^2ds-R_n\int_{\partial B_{Rn}^+} \Big|\frac{ \partial u}{\partial\nu}\Big|^{2}ds +\frac{R_n}2\l\int_{\partial B_{R_n}^+}u^2 ds\\
&\;\;\;-\frac{R_n}{p+1}\int_{\partial B_{R_n}^+}| u|^{p+1}ds -\frac{R_n}{q+1}\int_{\partial B_{R_n}\cap\Sigma_1}|u|^{q+1} ds,
\end{align*}
tends to $0$ as $n$ tends to $+\infty$.
 So we deduce from the Pohozaev identity \eqref {j} that
\begin{equation}\label{l}
\frac{N-2}{2}\int_{\R^N_+}|\nabla u|^{2}+\frac{N}{2}\l\int_{\R^N_+}u^2=\frac{N}{p+1}\int_{\R^N_+}|u|^{p+1}+\frac{N-1}{q+1}\int_{\Sigma_1}|u|^{q+1}dx'.
\end{equation}
On the other hand,  we multiply \eqref{11} by $u$ and integrate by parts, then we obtain
\begin{equation}\label{l1}
\int_{\R^N_+}|\nabla u|^{2}+\l\int_{\R^N_+}u^2=\int_{\R^N_+}|u|^{p+1}+\int_{\Sigma_1}|u|^{q+1}dx'.
\end{equation}
Combining \eqref {l} and \eqref {l1} we derive
\begin{equation}\label{l2}
\Big(\frac{N-2}{2}-\frac{N-1}{q+1}\Big)\int_{\R^N_+}|\nabla u|^2+\Big(\frac{N}{2}-\frac{N-1}{q+1}\Big)\l\int_{\R^N_+}u^2=\Big(\frac{N}{p+1}-\frac{N-1}{q+1}\Big)\int_{\R^N_+}|u|^{p+1}.
\end{equation}
Observe now since we assume that $p<\frac{N+2}{N-2}$ and $q>\frac{N}{N-2}$ then we have 

$$\frac{N-2}{2}-\frac{N-1}{q+1}>0  ,\quad \frac{N}{2}-\frac{N-1}{q+1}>0\quad \mbox{and} \;\;\;\frac{N}{p+1}-\frac{N-1}{q+1}>0.$$ 
So we deduce
$$\Big(\frac{N}{2}-\frac{N-1}{q+1}\Big)\l\int_{\R^N_+}u^2\leq\Big(\frac{N}{p+1}-\frac{N-1}{q+1}\Big)\int_{\R^N_+}|u|^{p+1}.$$
On the other hand by hypothesis $$|u|_{L^\infty(\R^N_+)}^{p-1}<\l(\frac{N}{2}-\frac{N-1}{q+1})(\frac{N}{p+1}-\frac{N-1}{q+1})^{-1},$$
then the last inequality implies
\begin{align*}
\Big(\frac{N}{2}-\frac{N-1}{q+1}\Big)\l\int_{\R^N_+}u^2\leq\Big(\frac{N}{p+1}-\frac{N-1}{q+1}\Big)|u|_{L^\infty(\R^N_+)}^{p-1}\int_{\R^N_+}u^2
<\Big(\frac{N}{2}-\frac{N-1}{q+1}\Big)\l\int_{\R^N_+}u^2.
\end{align*}
So, we reach a contradiction which completes the proof of Theorem \ref{t2}.\qed
\medskip

\subsection{\bf Proof of Theorem \ref{t3}.}
\medskip

We divide the proof in two parts.

\medskip
{\bf Step $1.$} If $p=\frac{N+2}{N-2}$. In view of the Pohozaev identity and as in the proof of Theorem \ref{t2}, we have
\begin{equation}\label{p0}
\frac{N-2}{2}\int_{\R^N_+}|\nabla u|^{2}+\frac{N}{2}\l\int_{\R^N_+}u^2=\frac{N-2}{2}\int_{\R^N_+}|u|^{\frac{2N}{N-2}}+\frac{N-1}{q+1}\int_{\Sigma_1}|u|^{q+1}.
\end{equation}
 On the other hand multiplying equation \eqref{11} or \eqref{12} by $u$ and integrating by parts yields
\begin{equation}\label{sq}
\int_{\R^N_+}|\nabla u|^{2}+\l\int_{\R^N_+}u^2=\int_{\R^N_+}|u|^{\frac{2N}{N-2}}+\int_{\Sigma_1}|u|^{q+1}.
\end{equation}
Combining \eqref{p0} and \eqref{sq} gives
\begin{equation}
\Big(\frac{N}{2}-\frac{N-2}{2}\Big)\l\int_{\R^N_+}u^2=\Big(\frac{N-1}{q+1}-\frac{N-2}{2}\Big)\int_{\Sigma_1}|u|^{q+1}.
\end{equation}
By assumption $2q-p-1\geq0$, then $\frac{N-1}{q+1}-\frac{N-2}{2}\leq0$. Hence we derive
$$\l\int_{\R^N_+}u^2=0.$$
Since $\l>0$ we must have $u\equiv 0$.
\medskip

{\bf Step $2.$} If $p>\frac{N+2}{N-2}.$ Substitute $R$ by $\tau R$ in \eqref{int11} (with $C_0> 0$) and scaling back, we obtain
\begin{equation*}
\int_{B_R^+}\bigg(|\nabla u^\tau|^2+|u^\tau|^{p+1}+\l\tau^2|u^\tau|^{2}\bigg)+\tau^{1-2\frac{q-1}{p-1}}\int_{\Sigma_1\cap B_R^+}|u^\tau|^{q+1}\leq C_0 \tau^{2\frac{p+1}{p-1}-N }+CR^{N-2\frac{p+1}{p-1}}.
\end{equation*}
As $2\frac{p+1}{p-1}-N<0$, then $(u^\tau)_{\tau \geq 1}$ is uniformly bounded in $H^1\cap L^{p+1}(B_R^+)$ for any $R>0$ and $(\tau^{\frac{p+1-2q}{(p-1)(q+1)}}u^\tau)_{\tau \geq 1}$ is uniformly bounded in $L^{q+1}(\Sigma_1\cap B_R^+)$. So, we can find
a sequence $(u^{\tau_j})$ which converges strongly in $L^{2}(B_R^+)$ and weakly in $H^1\cap L^{p+1}(B_R^+)$ to some function $u^\infty$ for every $R>0$, as $\tau_j\rightarrow +\infty$. Since $u^{\tau_j}$ satisfies the following PDE
\begin{eqnarray*}
  u^\tau=\lambda^{-1}\tau^{-2}\bigg(\D u^\tau + |u^\tau|^{p-1}u^\tau\bigg) \; &\text{ in \; $\R^{N}_{+}$},
\end{eqnarray*}
then taking limits in the sense of distributions, we get $u^\infty=0$, and therefore $(u^{\tau_j})$ converges strongly to $0$ in $L^{2}(B_R^+)$.
Invoking now inequality \eqref{int1} of Proposition \ref{pp} where we substitute $R$ by $\tau R$ we obtain
\begin{equation*}
\int_{B_{\tau R}^+}\bigg(|\nabla u|^2+|u|^{p+1}+\l u^{2}\bigg)+\int_{ \Sigma_1\cap B^+_{ \tau R }}|u|^{q+1}\leq C_0+ CR^{-2}\tau^{-2}\int_{B_{2\tau R}^+}u^2,\; \forall R>2R_0 \mbox{ and }\tau >1.
\end{equation*}
Thus, scaling back, we deduce
\begin{equation}\label{ew}
\int_{B_{ R}^+}\bigg(|\nabla u^\tau|^2+|u^\tau|^{p+1}+\l\tau^2|u^\tau|^{2}\bigg)+\tau^{1-2\frac{q-1}{p-1}}\int_{\Sigma_1\cap B_R^+}|u^\tau|^{q+1}\leq  C_0\tau^{2\frac{p+1}{p-1}-N }+CR^{-2}\int_{B_{2 R}^+}(u^\tau)^2.
\end{equation}
As a consequence, $(u^{\tau})$ converges strongly to $0$ in $H^1\cap L^{p+1}(B_R^+)$ and $(\tau u^{\tau})$ converges strongly to $0$ in $L^{2}(B_R^+)$ as $\tau\to \infty$.  Then from \eqref{ew}, we have
\begin{equation*}
\lim_{\tau\rightarrow +\infty}E_2(u,\tau)=0,
\end{equation*}
where $$E_2(u,\tau)=\int_{B_{1}^+}\bigg(\frac{1}{2}|\nabla u^\tau|^{2}+\frac{\lambda\tau^2}{2}|u^\tau|^{2}-\frac{1}{p+1}|u^\tau|^{p+1}\bigg)-\frac{\tau^{1-2\frac{q-1}{p-1}}}{q+1}\int_{\Sigma_1\cap B_1^+}|u^\tau|^{q+1}.$$
At this stage the rest of the proof is similar to the {\bf Case $3$} of Theorem \ref{t1}. Then we omit it. \qed

\medskip

\section{Acknowledgment}

The authors extend their appreciation to the Deanship of Scientific Research at King Khalid University, Abha, KSA for funding this work through Research Group under grant number (R.G.P-2 / 121/ 42).

\end{document}